\documentclass[15 pt]{article}
\usepackage{amsmath,amssymb,graphicx,bm}
\usepackage{graphicx} % Required for inserting images
\usepackage{epsfig}

\newtheorem{theorem}{Theorem}[section]
\newtheorem{lemma}[theorem]{Lemma}

\newtheorem{remark}[theorem]{Remark}

\title{Localized Orthogonal Decomposition Method with $H^1$ Interpolation for Multiscale Elliptic Problem}

\author{
Tao Yu
\footnote{School of Mathematics and Physics, Jinggangshan University, Ji'an, Jiangxi 343009, China. (yutao\_math@163.com)}
\and
Xingye Yue 
\footnote{Corresponding Author. Department of Mathematics, Soochow University, Suzhou, Jiangsu 215006, China. (xyyue@suda.edu.cn)}
}
\date{}

\begin{document}
\begin{sloppypar}
\maketitle

\begin{abstract}
This paper employs a localized orthogonal decomposition (LOD) method with $H^1$ interpolation for solving the multiscale elliptic problem. This method does not need any assumptions on scale separation. We give a priori error estimate for the proposed method. The theoretical results are conformed by various numerical experiments.
\end{abstract}

{\bf Keywords}\ \  {\small Localized Orthogonal Decomposition Method, $H^1$ Interpolation, Multiscale Elliptic Problem.}

\section{Introduction}
Consider the multiscale elliptic problem
\begin{eqnarray}\label{multiscale_elliptic_problem}
\left\{
\begin{array}{ll}\displaystyle
-\nabla\cdot(A^{\varepsilon}(x)\nabla u^{\varepsilon}(x))=f(x),\;\;&x\in\Omega,\\\displaystyle
u^{\varepsilon}(x)=0,&x\in\partial\Omega,
\end{array}
\right.
\end{eqnarray}
where $\Omega\subset\mathbb{R}^2$ (or $\mathbb{R}^3$) is a bounded convex polygonal domain with Lipschitz boundary $\partial \Omega$, $\varepsilon\ll 1$ is a positive parameter which signifies the multiscale nature of (\ref{multiscale_elliptic_problem}). For simplicity, we only consider $d=2$. All the results are still valid for high dimensional problems.

Direct numerical simulations such as finite element method (FEM) in order to capture fast oscillations on microscopic scale lead to problems of very large size which are prohibitively expensive. From the engineering point of view, the macroscopic properties and macroscopic behaviors of the materials are more important. Homogenization and numerical homogenization \cite{Bakhvalov-homogenization, Jikov-homogenization} have been used to simulate multiscale materials. These methods can calculate effective material properties to solve the macro scale problems successfully. However the limitation is based on the assumption of material, such as periodic micro-structures.

To overcome these difficulties, several multiscale computing strategies have been developed to solve the multiscale elliptic problem in heterogeneous structures on grids that are coarser than the scale of oscillations. Babu\v{s}ka and Osborn developed the so-called generalized multiscale finite element methods (GMSFEM) \cite{Babuska-GMSFEM1983, Babuska-GMSFEM1994,Babuska-GMSFEM2003} for problems with rough coefficients, including discontinuous coefficients as well as coefficients with multiple scales. Its main idea is to modify the finite element space in the framework of the finite element method. Hou et al. developed a multiscale finite element method (MsFEM) \cite{Efendiev-MsFEM, Hou-MsFEM} for the study of multi-dimensional problems. This is achieved by constructing multiscale finite element basis functions which adapt to the local property of the differential operator. E and Engquist developed the heterogeneous multiscale method (HMM) \cite{Abdulle-HMM, E-HMM1, E-HMM2}, which is a generally efficient methodology for problems with multiple scales. The method avoids computing the coefficients in the fine scale mesh as well as the coefficients of the homogenized equations. However, it depends on some powerful assumptions, such as periodicity and scale separation. There are some other methods to solve multiscale problems with strongly varying coefficients. For example, the variational multiscale method (VMS) \cite{Hughes-VMS1, Hughes-VMS2} and residual-free bubbles method \cite{Franca-RFB} in which the test function space was decomposed into the sum of coarse and fine scale components.

In this paper, we study the Localized Orthogonal Decomposition (LOD) method with $H^1$ interpolation for elliptic problems with highly varying coefficients. The LOD method was first introduced for elliptic multiscale problems \cite{Malqvist-LOD-elliptic}, which is suitable for the general framework of the Variational Multiscale Method (VMS)\cite{Hughes-VMS1, Hughes-VMS2}. The method constructs local generalized finite element basis which is exponentially decayed by using the modified Cl\'{e}ment interpolation \cite{Carstensen-interpolation}. The analysis does not depend on the regularity of the solution or the scale separation in the coefficient. The modified Cl\'{e}ment interpolation is essentially a $L^2$ projection, but $H^1$ interpolation is more natural for elliptic problems. The method is also applied to parabolic problem \cite{Malqvist-LOD-parabolic}, semi-linear problem \cite{Henning-LOD-semilinear}, wave propagation \cite{Abdulle-LOD-wave}, elliptic optimal control problem \cite{Brenner-LOD-OCP}, and so on. Recently, Peterseim reviewed the VMS methods of Linear multiscale partial differential equations \cite{Peterseim-VMS-decay} and gave a new simple proof for the exponential decay property. We will use $H^1$ interpolation to apply the technique. By using the classical finite element methods, we prove convergence of optimal order in $H^1$ norm. Numerical examples are given to complete the analysis, which supports our theoretical findings.

The outline of this paper is as follows. In Section 2, we describe the model problem and present the localized orthogonal decomposition method with $H^1$ interpolation. The localization technique is proposed in Section 3. In section 4, some numerical experiments are given.

\section{LOD Method with $H^1$ Interpolation for Multiscale Elliptic Problem}
\setcounter{equation}{0}
In this section, we first construct a $H^1$ interpolation similar to the modified Cl\'{e}ment interpolation in \cite{Carstensen-interpolation}. Then, we introduce the localized orthogonal decomposition method with this $H^1$ interpolation for the multiscale elliptic problem (\ref{multiscale_elliptic_problem}).

Assume that $f\in L^2(\Omega)$, and $A^{\varepsilon}(x)=(A^{\varepsilon}_{ij}(x))$ is a positive definite, bounded tensor satisfying
\begin{eqnarray}
\label{A-property}
{\alpha|\xi|^2\leq
A^{\varepsilon}_{ij}\xi_i\xi_j\leq \beta|\xi|^2\;,\;\;\forall
\xi\in \mathbb{R}^d,}
\end{eqnarray}
for some positive constants $\alpha$ and $\beta$. Multiply equation (\ref{multiscale_elliptic_problem}) by test function $v\in V \overset{\vartriangle}{=} H^1_0(\Omega)$, by using the Green's formula, the standard weak formulation is to find $u^{\varepsilon}\in V$ such that
\begin{eqnarray}
\label{weak form}
a(u^{\varepsilon},v)\overset{\vartriangle}{=}\int_{\Omega}A^{\varepsilon}(x)\nabla u^{\varepsilon}(x)\cdot\nabla v(x)dx =\int_{\Omega}f(x)v(x)dx\overset{\vartriangle}{=}(f,v),\;\;\forall v \in V.
\end{eqnarray}
The bilinear form $a(\cdot,\cdot)$ is symmetric, continuous and coercive, directly from property (\ref{A-property}) and Poincar\'{e} inequality. According to the Lax-Milgram Theorem, the weak formulation (\ref{weak form}) has a unique weak solution $u^{\varepsilon}\in V$, satisfying
\begin{eqnarray*}
||u^{\varepsilon}||_V\leq C ||f||_{0, \Omega}.
\end{eqnarray*}
Let $\mathcal{T}_h$ be a standard regular coarse triangulation of the domain $\Omega$ with coarse mesh size $h_K\overset{\vartriangle}{=}diam(K)$ for all $K\in \mathcal{T}_h$. Denote by $h=\underset{K\in \mathcal{T}_h}{\max}h_K$ and $\mathcal{N}$ the set of all interior nodal points, $\mathcal{N}(K)$ the vertices of $K\in \mathcal{T}_h$. The classical conforming $P_1$ finite element space is given by
\begin{eqnarray*}
V_h=\{u_h\in H^1_0(\Omega)|\;u_h|_{K}\in P_1(K),\;\forall K\in \mathcal{T}_h\},
\end{eqnarray*}
where $P_1(K)$ is the space of polynomials with total degree at most one on element $K$. Then the finite element method reads: find $u_h\in V_h$ such that
\begin{eqnarray}
a(u_h,v_h)=(f,v_h),\;\forall v_h\in V_h.
\end{eqnarray}
For every interior nodal point $a\in \mathcal{N}$, let $\lambda_a\in V_h$ be the corresponding nodal basis function, i.e. the continuous piecewise linear finite element function which equals one at $a$ and vanishes at all other vertices. In order to construct the modified nodal basis, we present an important $H^1$ interpolation $\mathcal{J}_h: V \rightarrow V_h$. For every $v\in V$, define
\begin{eqnarray}
\label{H1_interpolation}
\mathcal{J}_h v=\sum_{a\in \mathcal{N}} \tilde{v}(a)\lambda_a,
\end{eqnarray}
where
\begin{eqnarray}\label{tildav}
\tilde{v}(a)=\frac{\int_{\Omega}v \lambda_a dx + h^2 \displaystyle\sum_{K\in \mathcal{T}_h}\int_{K}\nabla v \cdot \nabla\lambda_a dx}{\int_{\Omega} \lambda_a dx + h^2 \displaystyle\sum_{K\in \mathcal{T}_h}\int_{K}|\nabla\lambda_a| dx}.
\end{eqnarray}
For the above special definition, the following bounded property can be easily derived.
\begin{lemma}\label{bounded property}
There exists a positive constant $C$, which is independent on the mesh size $h$, such that for every interior nodal point $a\in K\in \mathcal{T}_h$ and $v\in V$, it holds that
\begin{eqnarray}
||\tilde{v}(a)||_{L^2(K)}\leq C\left(||v||_{L^2(w_a)}+h||\nabla v||_{L^2(w_a)}\right),
\end{eqnarray}
where $w_a \overset{\vartriangle}{=} \bigcup \{ K\in \mathcal{T}_h|a \in K\}$.  
\end{lemma}

\textbf{Proof.} Using the definition (\ref{tildav}) of $\tilde{v}(a)$, it shows that
\begin{eqnarray}
||\tilde{v}(a)||_{L^2(K)}\!\!\!\!&=&\!\!\!\!|K|^{\frac{1}{2}}|\tilde{v}(a)| =|K|^{\frac{1}{2}}\left|\frac{\int_{\Omega} v \lambda_a dx + h^2 \displaystyle\sum_{K\in \mathcal{T}_h}\int_{K} \nabla v \cdot \nabla\lambda_a dx}{\int_{\Omega} \lambda_a dx + h^2 \displaystyle\sum_{K\in \mathcal{T}_h}\int_{K}|\nabla\lambda_a| dx}\right|\nonumber\\
\!\!\!\!&\leq&\!\!\!\!|K|^{\frac{1}{2}}\left(\int_{\omega_a} \lambda_a dx\right)^{-1}\left(\left|\int_{\omega_a} v \lambda_a dx\right|+\left|h^2 \displaystyle\sum_{K\in \omega_a}\int_{K}\nabla v \cdot \nabla\lambda_a dx\right| \right)\nonumber\\
\!\!\!\!&\leq&\!\!\!\!|K|^{\frac{1}{2}}\left(\frac{1}{3}|\omega_a|\right)^{-1} \Bigg(||\lambda_a||_{L^2(\omega_a)}||v||_{L^2(\omega_a)}\nonumber\\
\!\!\!\!&&\!\!\!\!\qquad\qquad\qquad\qquad+h^2\left(\sum_{K\in \omega_a}||\nabla\lambda_a||^2_{L^2(K)}\right)^{\frac{1}{2}}||\nabla v||_{L^2(\omega_a)}\Bigg)\nonumber\\
\!\!\!\!&\leq&\!\!\!\!C |K|^{\frac{1}{2}}|\omega_a|^{-1}\left(|\omega_a|^{\frac{1}{2}}||v||_{L^2(\omega_a)}+h^2 h^{-1}|\omega_a|^{\frac{1}{2}}||\nabla v||_{L^2(\omega_a)}\right)\nonumber\\
\!\!\!\!&=&\!\!\!\!C|K|^{\frac{1}{2}}|\omega_a|^{-\frac{1}{2}} \left(||v||_{L^2(\omega_a)}+h||\nabla v||_{L^2(\omega_a)} \right)\nonumber\\
\!\!\!\!&\leq&\!\!\!\!C \left(||v||_{L^2(\omega_a)}+h||\nabla v||_{L^2(\omega_a)} \right),\nonumber
\end{eqnarray}
directly by H\"{o}lder inequality and simple derivation.\hfill $\blacksquare$\\
Next, we can obtain the quasi-local stability and approximation properties for $H^1$ interpolation $\mathcal{J}_h$ as follow.
\begin{lemma}
\label{interpolation property}
There exists a positive constant $C_{\mathcal{J}_h}$, which is independent on the mesh size $h$, such that for all $v\in V$ and $K\in \mathcal{T}_h$, it holds that
\begin{eqnarray}\label{local_stability}
||\mathcal{J}_h v||_{L^2(K)}\leq C_{\mathcal{J}_h} ||v||_{H^1(w_K)}
\end{eqnarray}
and
\begin{eqnarray}\label{approximation_property}
||v-\mathcal{J}_h v||_{L^2(K)}\leq C_{\mathcal{J}_h} h||v||_{H^1(w_K)},
\end{eqnarray}
where $w_K \overset{\vartriangle}{=} \bigcup \{ T\in \mathcal{T}_h|T\bigcap K\neq \emptyset\}$.
\end{lemma}
\textbf{Proof.} Since $\displaystyle\sum_{a\in \mathcal{N}(K)}\lambda_a=1$ on $K\in \mathcal{T}_h$, we conclude from the definition (\ref{H1_interpolation}) of $H^1$ interpolation $\mathcal{J}_h$ that
\begin{eqnarray}
||v-\mathcal{J}_h v||_{L^2(K)}\!\!\!\!&=&\!\!\!\!||\sum_{a\in \mathcal{N}(K)}\lambda_a\left(v-\tilde{v}(a)\right)||_{L^2(K)}\leq\sum_{a\in \mathcal{N}(K)}||v-\tilde{v}(a)||_{L^2(K)}\nonumber\\
\!\!\!\!&\leq&\!\!\!\!\sum_{a\in \mathcal{N}(K)}C h||v||_{H^1(w_K)}=C_{\mathcal{J}_h} h||v||_{H^1(w_K)},
\end{eqnarray}
in which we need to prove the inequality
\begin{eqnarray}\label{v-tildev0}
||v-\tilde{v}(a)||_{L^2(K)}\leq C h||v||_{H^1(w_K)}.
\end{eqnarray}
In fact,
\begin{eqnarray}\label{v-tildev}
||v-\tilde{v}(a)||_{L^2(K)}\leq ||v-\bar{v}||_{L^2(K)}+||\bar{v}-\tilde{v}(a)||_{L^2(K)},
\end{eqnarray}
where $\bar{v}=\frac{1}{|K|}\int_{K}vdx$ is the mean-value of $v$ on the element $K$.
From a standard Poincare-Friedrichs inequality, we know that
\begin{eqnarray}\label{v-barv}
||v-\bar{v}||_{L^2(K)}\leq C diam(K)||\nabla v||_{L^2(K)}\leq C h_K||\nabla v||_{L^2(w_K)}.
\end{eqnarray}
For the second term of (\ref{v-tildev}), denote by $\pi_a v = \tilde{v}(a)$, then
\begin{eqnarray}\label{pi_a_bar_v}
\pi_a\bar{v}\!\!\!\!&=&\!\!\!\!\frac{\int_{\Omega}\bar{v} \lambda_a dx + h^2 \displaystyle\sum_{K\in \mathcal{T}_h}\int_{K}\nabla \bar{v} \cdot \nabla\lambda_a dx}{\int_{\Omega} \lambda_a dx + h^2 \displaystyle\sum_{K\in \mathcal{T}_h} \int_{K}|\nabla\lambda_a| dx}=\frac{\bar{v} \int_{\Omega}\lambda_a dx}{\int_{\Omega} \lambda_a dx + h^2 \displaystyle\sum_{K\in \mathcal{T}_h}\int_{K}|\nabla\lambda_a| dx}\nonumber\\
\!\!\!\!&=&\!\!\!\!\frac{\bar{v} \left(\int_{\Omega}\lambda_a dx+h^2 \displaystyle\sum_{K\in \mathcal{T}_h}\int_{K}|\nabla\lambda_a| dx\right) - \bar{v} h^2 \displaystyle\sum_{K\in \mathcal{T}_h}\int_{K}|\nabla\lambda_a| dx}{\int_{\Omega} \lambda_a dx + h^2 \displaystyle\sum_{K\in \mathcal{T}_h}\int_{K}|\nabla\lambda_a| dx}\nonumber\\
\!\!\!\!&=&\!\!\!\!\bar{v}-\left(\frac{h^2 \displaystyle\sum_{K\in \mathcal{T}_h}\int_{K}|\nabla\lambda_a| dx}{\int_{\Omega} \lambda_a dx + h^2\displaystyle\sum_{K\in \mathcal{T}_h} \int_{K}|\nabla\lambda_a| dx}\right)\bar{v}.\nonumber
\end{eqnarray}
Thus
\begin{eqnarray}\label{barv-tildev}
\!\!\!\!&&\!\!\!\!||\bar{v}-\tilde{v}(a)||_{L^2(K)}=||\bar{v}-\pi_a v||_{L^2(K)}\nonumber\\
\!\!\!\!&=&\!\!\!\!||\pi_a\bar{v}+\left(\frac{h^2 \displaystyle\sum_{K\in \mathcal{T}_h}\int_{K}|\nabla\lambda_a| dx}{\int_{\Omega} \lambda_a dx + h^2 \displaystyle\sum_{K\in \mathcal{T}_h}\int_{K}|\nabla\lambda_a| dx}\right)\bar{v}-\pi_a v||_{L^2(K)}\nonumber\\
\!\!\!\!&\leq&\!\!\!\!||\pi_a\bar{v}-\pi_a v||_{L^2(K)}+||\left(\frac{h^2 \displaystyle\sum_{K\in \mathcal{T}_h}\int_{K}|\nabla\lambda_a| dx}{\int_{\Omega} \lambda_a dx + h^2 \displaystyle\sum_{K\in \mathcal{T}_h}\int_{K}|\nabla\lambda_a| dx}\right)\bar{v}||_{L^2(K)}.
\end{eqnarray}
The first term of (\ref{barv-tildev}) can be deduced by
\begin{eqnarray}\label{pibarv-piv}
||\pi_a\bar{v}-\pi_a v||_{L^2(K)}=||\pi_a(v-\bar{v}) ||_{L^2(K)}
\leq C \left(||v-\bar{v}||_{L^2(\omega_a)}+h||\nabla(v-\bar{v})||_{L^2(\omega_a)} \right)\nonumber\\
\leq C \left(diam(\omega_a)||\nabla v||_{L^2(\omega_a)}+h||\nabla v||_{L^2(\omega_a)} \right)\leq C h ||\nabla v||_{L^2(\omega_a)},
\end{eqnarray}
in which Lemma \ref{bounded property} and (\ref{v-barv}) are used. Using the definition of $\bar{v}$ and the H\"{o}lder inequality, it easily see that
\begin{eqnarray}\label{residual term}
\!\!\!\!&&\!\!\!\!||\left(\frac{h^2 \displaystyle\sum_{K\in \mathcal{T}_h}\int_{K}|\nabla\lambda_a| dx}{\int_{\Omega} \lambda_a dx + h^2 \displaystyle\sum_{K\in \mathcal{T}_h}\int_{K}|\nabla\lambda_a| dx}\right)\bar{v}||_{L^2(K)}\nonumber\\
\!\!\!\!&=&\!\!\!\!||\left(\frac{h^2 \displaystyle\sum_{K\in \omega_a}\int_{K}|\nabla\lambda_a| dx}{\int_{\omega_a} \lambda_a dx + h^2 \displaystyle\sum_{K\in \omega_a}\int_{K}|\nabla\lambda_a| dx}\right)\bar{v}||_{L^2(K)}\nonumber\\
\!\!\!\!&\leq&\!\!\!\! \frac{h^2 \displaystyle\sum_{K\in \omega_a}\int_{K}|\nabla\lambda_a| dx}{\int_{\omega_a} \lambda_a dx}|\bar{v}||K|^{\frac{1}{2}}
\leq C\frac{h^2 |\omega_a|h^{-1}}{|\omega_a|}|K|^{-1}|K|^{\frac{1}{2}}||v||_{L^2(K)}|K|^{\frac{1}{2}}\nonumber\\
\!\!\!\!&\leq&\!\!\!\! Ch||v||_{L^2(K)}.
\end{eqnarray}
Inequalities (\ref{v-tildev})-(\ref{residual term}) immediately yield the desired result (\ref{v-tildev0}). Obviously from Lemma \ref{bounded property},
\begin{eqnarray}\label{pi_a_v}
||\tilde{v}(a)||_{L^2(K)}\leq C\left(||v||_{L^2(\omega_a)}+h||\nabla v||_{L^2(\omega_a)}\right)\leq C||v||_{H^1(\omega_a)}.
\end{eqnarray}
Therefore
\begin{eqnarray}\label{pi_a_v_1}
||\mathcal{J}_h v||_{L^2(K)}\!\!\!\!&=&\!\!\!\!||\sum_{a\in \mathcal{N}} \tilde{v}(a)\lambda_a||_{L^2(K)} \leq \sum_{a\in \mathcal{N}}||\tilde{v}(a)||_{L^2(K)}\nonumber\\
\!\!\!\!&\leq&\!\!\!\! \sum_{a\in \mathcal{N}}C||v||_{H^1(\omega_a)}\leq C_{\mathcal{J}_h} ||v||_{H^1(w_K)}
\end{eqnarray}
is the desired result (\ref{local_stability}).\hfill $\blacksquare$

The kernel of $\mathcal{J}_h$,
\begin{eqnarray}
\label{kernel}
V^f\overset{\vartriangle}{=}\{v\in V|\mathcal{J}_h\;v=0\}
\end{eqnarray}
represents the microscopic features of $V$, which are not captured by $V_h$. Given any nodal basis $\lambda_a$ of $V_h$, define fine-scale corrector $\phi_a=\mathfrak{F}\lambda_a\in V^f$, satisfying the corrector problem
\begin{eqnarray}
\label{corrector problem}
a(\phi_a,\omega)=a(\lambda_a,\omega),\;\forall\;\omega\in V^f.
\end{eqnarray}
The finescale projection operator $\mathfrak{F}: V_h \rightarrow V^f$ leads to an orthogonal splitting
\begin{eqnarray}\label{orthogonal splitting}
V=V_h^{ms}\oplus V^f,
\end{eqnarray}
with respect to the bilinear form $a(\cdot, \cdot)$, where $V_h^{ms}=V_h-\mathfrak{F} V_h$ is called modified coarse space, containing fine scale information.
That is to say, for any $u\in V$, there exist $u_h^{ms}\in V_h^{ms}$ and $u^f\in V^f$, such that
\begin{eqnarray}
u=u_h^{ms}+u^f,\;\;a(u_h^{ms},u^f)=0.
\end{eqnarray}
The basis of $V_h^{ms}$ is then given by the modified nodal basis $\lambda_a-\phi_a$, i.e.,
\begin{eqnarray}
V_h^{ms}=span\{\lambda_a-\phi_a|a\in \mathcal{N}\}.
\end{eqnarray}
Hence, the corresponding Galerkin approximation solution $u_h^{ms}\in V_h^{ms}$ satisfies
\begin{eqnarray}
\label{Galerkin approximation}
a(u_h^{ms},v)=(f,v),\;\;\forall v\in V_h^{ms}.
\end{eqnarray}
The error $u^{\varepsilon}-u_h^{ms}$ of the above method (\ref{Galerkin approximation}) is analyzed by the following Theorem.
\begin{theorem}\label{macro_error}
 Let $u^{\varepsilon}\in V$ and $u_h^{ms}\in V_h^{ms}$ be the solutions of problems (\ref{weak form}) and (\ref{Galerkin approximation}), respectively. Then it holds that
\begin{eqnarray}
||u^{\varepsilon}-u_h^{ms}||_{V}\leq C h ||f||_{L^2(\Omega)},
\end{eqnarray}
where the positive constant $C$ is not dependent on the mesh size $h$.
\end{theorem}
\textbf{Proof.} Let $v=v_h^{ms}\in V_h^{ms}$ in (\ref{weak form}) and (\ref{Galerkin approximation}), respectively. The orthogonality results
\begin{eqnarray}\label{orthogonality}
a(u^{\varepsilon}-u_h^{ms},v_h^{ms})=0,
\end{eqnarray}
can be easily obtained by subtracting these two problems.
Denote $u^{\varepsilon}-u_h^{ms}=u^f$, then $\mathcal{J}_h u^f=0$.
From the application of (\ref{A-property}) and the Poinc\'{a}re inequality, it yield
\begin{eqnarray}\label{a1}
a(u^f,u^f)=\int_{\Omega}A^{\varepsilon} \nabla u^f\cdot \nabla u^f dx \geq \alpha \int_{\Omega}\left(\nabla u^f\right)^2 dx\geq\frac{\alpha}{\tilde{C}} ||u^f||^2_{V}.
\end{eqnarray}
On the other hand,
\begin{eqnarray}\label{a2}
a(u^f,u^f)\!\!\!\!&=&\!\!\!\! a(u^{\varepsilon},u^f)-a(u_h^{ms},u^f)=a(u^{\varepsilon},u^f)=(f,u^f)\nonumber\\
\!\!\!\!&=&\!\!\!\!\sum_{K}\int_{K}fu^fdx\leq \sum_{K}||f||_{L^2(K)}||u^f||_{L^2(K)}\nonumber\\
\!\!\!\!&=&\!\!\!\!\sum_{K}||f||_{L^2(K)}||u^f-\mathcal{J}_h u^f||_{L^2(K)}\nonumber\\
\!\!\!\!&\leq&\!\!\!\!\sum_{K}||f||_{L^2(K)}C_{\mathcal{J}_h}h_K ||u^f||_{V,\;\omega_K}\nonumber\\
\!\!\!\!&\leq&\!\!\!\!\sum_{K}\left(\frac{C_{\mathcal{J}_h}^2 h_K^2}{2\delta}||f||^2_{L^2(K)} +\frac{\delta}{2}||u^f||^2_{V,\;\omega_K}\right)\nonumber\\
\!\!\!\!&\leq&\!\!\!\!\frac{C_{\mathcal{J}_h}^2 h^2}{2\delta}||f||^2_{L^2(\Omega)} +\frac{C\delta}{2}||u^f||^2_{V},
\end{eqnarray}
by using Lemma \ref{interpolation property}, (\ref{orthogonality}) and Young's inequality. Thus, together with (\ref{a1}) and (\ref{a2}),
\begin{eqnarray}
\frac{\alpha}{\tilde{C}} ||u^f||^2_{V}\leq \frac{C_{\mathcal{J}_h}^2 h^2}{2\delta}||f||^2_{L^2(\Omega)} +\frac{C\delta}{2}||u^f||^2_{V}.
\end{eqnarray}
The small enough choice for $\delta$  concludes the proof.\hfill $\blacksquare$

\section{Localization}
In this section, we first show that the fine-scale corrector $\phi_a$ has exponential decay properties, by using the technique in \cite{Peterseim-VMS-decay}. Therefore, simple truncation of the corrector problems to local patches of coarse elements yields localized basis functions with good approximation properties.

Let $\eta\in W^{1,\infty}(\Omega)$ be some cut-off function, satisfying
\begin{eqnarray}
\label{cut_off_function}
\left\{
\begin{array}{ll}\displaystyle
\eta=1,& x\in B_r(a),\\
\eta=0,& x\in\Omega\setminus B_R(a),
\end{array}
\right.~~and~~||\nabla\eta||_{L^{\infty}(\Omega)}\leq C_{\eta}h^{-1},
\end{eqnarray}
with $R>2h$, $r=R-h>h$ and constant $C_{\eta}$ independent of $h$, where $B_r(a)$ and $B_R(a)$ are  balls of radius $r, R>0$ centred at $a$. Multiplying by $\eta$ and $(1-\mathcal{J}_h)$, function in the kernel of $\mathcal{J}_h$ has quasi-local stable property in the sense that
\begin{eqnarray}\label{quasi_local_stable}
||(1-\mathcal{J}_h)(\eta\omega)||_{H^{1}, B_R(a)\setminus B_r(a)}\leq C_{\eta,\mathcal{J}_h}||\omega||_{H^{1}, B_{R'}(a)\setminus B_{r'}(a)},\;\forall\;\omega\in Ker(\mathcal{J}_h)
\end{eqnarray}
holds with $r'=r-mh$ and $R'=R+mh$, in which constant $C_{\eta,\mathcal{J}_h}>0$ and $m\in N_0$ independent of $h$ and $a\in \mathcal{N}$.
This result is possible because $\mathcal{J}_h$ enjoys quasi-local stability and approximation properties in Lemma \ref{interpolation property}.
Similarly as \cite{Peterseim-VMS-decay}, we can see that the finescale corrector $\phi_a$ decays exponentially.
\begin{lemma}\label{Exponential decay} There exist constants $c>0$ and $C>0$, independent of $h$ and $R$, such that
\begin{eqnarray}\label{Exponential decay1}
||\phi_a||_{H^{1},\;\Omega \setminus B_R(a)}\leq C\;\exp\left(-c\frac{R}{h}\right)||\phi_a||_{H^{1},\;\Omega},
\end{eqnarray}
where $B_R(a)$ is a ball of radius $R>0$ centred at $a$.
\end{lemma}
\textbf{Proof.} Because the finescale corrector $\phi_a \in V^f$, i.e. $\mathcal{J}_h \phi_a = 0$, it satisfies
\begin{eqnarray*}
||\phi_a||_{H^{1},\;\Omega\setminus B_R(a)}  \!\!\!\!&=&\!\!\!\!||(1-\mathcal{J}_h)\phi_a||_{H^{1},\;\Omega\setminus B_R(a)} \\
\!\!\!\!&=&\!\!\!\!||(1-\mathcal{J}_h)(1-\eta) \phi_a||_{H^{1},\;\Omega\setminus B_R(a)}\\
\!\!\!\!&\leq&\!\!\!\!||(1-\mathcal{J}_h)(1-\eta) \phi_a||_{H^{1},\;\Omega}\\
\!\!\!\!&\leq&\!\!\!\!\beta^{-1} a(\omega, (1-\mathcal{J}_h)(1-\eta) \phi_a)\\
\!\!\!\!&=&\!\!\!\!\beta^{-1} \left(a(\omega,\phi_a)-a(\omega,(1-\mathcal{J}_h)\eta \phi_a) \right),
\end{eqnarray*}
where $\omega=\frac{(1-\mathcal{J}_h)(1-\eta) \phi_a }{||(1-\mathcal{J}_h)(1-\eta) \phi_a||_{V, \Omega}} \in V^f$ satisfies $||w||_{H^{1},\;\Omega}=1$ and
\begin{eqnarray}\label{suppt_w}
supp\;\omega = supp\left( (1-\mathcal{J}_h)(1-\eta) \phi_a\right) \subset \Omega\setminus B_r(a).
\end{eqnarray}
It's obviously that the supports of $\lambda_a$ and $\omega$ have no overlap, then
\begin{eqnarray*}
a(\omega,\phi_a)=a(\phi_a,\omega)=a(\lambda_a,\omega)=0,
\end{eqnarray*}
by using the symmetry property of $a(\cdot,\cdot)$ and (\ref{corrector problem}). Thus
\begin{eqnarray}\label{phi_a_BR}
||\phi_a||_{H^{1},\;\Omega\setminus B_R(a)}\leq \beta^{-1}\left|a(\omega,(1-\mathcal{J}_h)\eta \phi_a)\right|.
\end{eqnarray}
Directly from (\ref{cut_off_function}), (\ref{quasi_local_stable}), (\ref{suppt_w}), (\ref{phi_a_BR}) and Lemma \ref{interpolation property},
\begin{eqnarray*}
||\phi_a||_{H^{1},\;\Omega\setminus B_{R'}(a)}
\!\!\!\!&\leq&\!\!\!\!||\phi_a||_{H^{1},\;\Omega\setminus B_R(a)}\leq\beta^{-1} |a(\omega,(1-\mathcal{J}_h)\eta \phi_a)|\\
\!\!\!\!&\leq&\!\!\!\!\alpha \beta^{-1} ||(1-\mathcal{J}_h)\eta \phi_a|| _{H^{1},\;\Omega\setminus B_r(a)}\\ \!\!\!\!&=&\!\!\!\!\alpha\beta^{-1} ||(1-\mathcal{J}_h)\eta\phi_a|| _{H^{1},\;B_{R}(a)\setminus B_{r}(a)}\\
\!\!\!\!&\leq&\!\!\!\!\alpha\beta^{-1} C_{\eta, \mathcal{J}_h}||\phi_a||_{H^{1},\;B_{R'}(a)\setminus B_{r'}(a)}\\
\!\!\!\!&=&\!\!\!\!\alpha\beta^{-1}C_{\eta, \mathcal{J}_h}\left(||\phi_a||^2_{H^{1},\;\Omega\setminus B_{r'}(a)}- ||\phi_a||^2_{H^{1},\;\Omega\setminus B_{R'}(a)}\right)^{\frac{1}{2}}.
\end{eqnarray*}
Thus
\begin{eqnarray*}
||\phi_a||^2_{H^{1},\;\Omega\setminus B_{R'}(a)}\leq \frac{C'}{1+C'}||\phi_a||^2_{H^{1},\;\Omega\setminus B_{r'}(a)},
\end{eqnarray*}
where $C'=\left(\alpha\beta^{-1}C_{\eta, \mathcal{J}_h}\right)^2$. The iterative application of this estimate with $\frac{R}{(2m+1)h}$, and relabelling $R'$ to $R$, have
\begin{eqnarray*}
 ||\phi_a||^2_{H^{1},\;\Omega\setminus B_{R}(a)}\leq\left( \frac{C'}{1+C'}\right)^{\frac{R}{(2m+1)h}}||\phi_a||^2_{H^{1},\;\Omega} = \exp(-2c\frac{R}{h})||\phi_a||^2_{H^{1},\;\Omega}
\end{eqnarray*}
where $c=\left|ln(\frac{C'}{1+C'})\right|\frac{1}{2(2m+1)}>0$.\hfill $\blacksquare$

Hence, simple truncation of the corrector problems to local patches of coarse elements yields localized basis functions with good approximation properties. For positive integer $l$, define $l$-th order element patch $\Omega_{K,l}$ about $K\in \mathcal{T}_h$ by
\begin{eqnarray}
&&\Omega_{K,1}\overset{\vartriangle}{=}int\left(\bigcup \{T\in\mathcal{T}_h|T\cap K\neq\emptyset\}\right),\\
&&\Omega_{K,l}\overset{\vartriangle}{=}int\left(\bigcup \{T\in\mathcal{T}_h|T\cap \bar{\Omega}_{K,l-1}\neq\emptyset\}\right).
\end{eqnarray}
Then for all $u,\; v\in V$ and $\omega\subset\Omega$, the localized bilinear form reads
\begin{eqnarray}
a_{\omega}(u,v)\overset{\vartriangle}{=}\int_{\omega}A^{\varepsilon}(x)\nabla u(x)\cdot\nabla v(x) dx.
\end{eqnarray}
Given any nodal basis function $\lambda_a\in V_h$, let $\phi_{a,l,K}\in V^f \cap H_0^1(\Omega_{K,l})$ solve the subscale corrector problem
\begin{eqnarray}\label{localised corrector problem}
a_{\Omega_{K,l}}(\phi_{a,l,K},w)=a_{K}(\lambda_a,w),\;\; \forall w\in V^f \cap H_0^1(\Omega_{K,l}).
\end{eqnarray}
Let $\phi_{a,l}=\underset{K\in \mathcal{T}_h,\;a\in K}{\sum}\phi_{a,l,K}$ and the localized corrector space
\begin{eqnarray}
V_{h,l}^{ms}\overset{\vartriangle}{=}span\{\lambda_a-\phi_{a,l}|a\in \mathcal{N}\}.
\end{eqnarray}
The localised multiscale Galerkin FEM seeks $u_{h,l}^{ms}\in V_{h,l}^{ms}$ such that
\begin{eqnarray}\label{localised method}
a(u_{h,l}^{ms},w_{h,l}^{ms})=(f,w_{h,l}^{ms}),\;\;\forall w_{h,l}^{ms}\in V_{h,l}^{ms}.
\end{eqnarray}
\begin{remark}
The trial function space $V_h^{ms}$ and $V_{h,l}^{ms}$ in (\ref{Galerkin approximation}) and (\ref{localised method}) can be reduced to $V_h$, directly from the corrector problems (\ref{corrector problem}) and (\ref{localised corrector problem}).
\end{remark}
The error estimate for the localized method (\ref{localised method}) is analysed by the following Theorem.
\begin{theorem}\label{Error estimate} There exists positive constant $C$, independent of $h$, such that
\begin{eqnarray}
||u^{\varepsilon}-u_{h,l}^{ms}||_{V}\leq Ch||f||_V,
\end{eqnarray}
where $u^{\varepsilon}\in V$ and $u_{h,l}^{ms}\in V_{h,l}^{ms}$ be the solutions of problems (\ref{weak form}) and (\ref{localised method}), respectively.
\end{theorem}
\textbf{Proof.} The error can be decomposed into two parts
\begin{eqnarray}
||u^{\varepsilon}-u_{h,l}^{ms}||_{V}\leq||u^{\varepsilon}-u_{h}^{ms}||_{V}+ ||u_{h}^{ms}-u_{h,l}^{ms}||_{V}.
\end{eqnarray}
The first term is called as the macro error
\begin{eqnarray}
||u^{\varepsilon}-u_{h}^{ms}||_{V}\leq Ch||f||_V,
\end{eqnarray}
directly from Theorem \ref{macro_error}.
It remains to estimate the second term $||u_{h}^{ms}-u_{h,l}^{ms}||_{V}$. Let $Z_h^{ms}\in V_{h}^{ms}$ and $||Z_h^{ms}||_V=1$, then
\begin{eqnarray}
\alpha ||u_{h}^{ms}-u_{h,l}^{ms}||_{V}\!\!\!\!&\leq&\!\!\!\! a(u_{h}^{ms}-u_{h,l}^{ms},Z_h^{ms})= a(u_{h}^{ms}-u_{h,l}^{ms},Z_h^{ms}-Z_{h,l}^{ms})\nonumber\\
\!\!\!\!&=&\!\!\!\!a(u_{h}^{ms}-u^{\varepsilon},Z_h^{ms}-Z_{h,l}^{ms})+ a(u^{\varepsilon}-u_{h,l}^{ms},Z_h^{ms}-Z_{h,l}^{ms})\nonumber\\
\!\!\!\!&=&\!\!\!\!a(u^{\varepsilon}-u_{h,l}^{ms},Z_h^{ms}-Z_{h,l}^{ms})\nonumber\\
\!\!\!\!&\leq&\!\!\!\!\beta ||u^{\varepsilon}-u_{h,l}^{ms}||_V||Z_h^{ms}-Z_{h,l}^{ms}||_V
\end{eqnarray}
The exponential decay property (\ref{Exponential decay1}) allows one to choose $Z_{h,l}^{ms}$ in such a way that
\begin{eqnarray}
||Z_h^{ms}-Z_{h,l}^{ms}||_V\leq C\;\exp\left(-c l\right)||Z_h^{ms}||_V =C\;\exp\left(-c l\right).
\end{eqnarray}
Hence
\begin{eqnarray}
||u_{h}^{ms}-u_{h,l}^{ms}||_{V}\leq \frac{\beta}{\alpha}C\;\exp\left(-cl\right) ||u^{\varepsilon}-u_{h,l}^{ms}||_V.
\end{eqnarray}
Choose $l$ big enough, such that
\begin{eqnarray}
\frac{\beta}{\alpha}C\;\exp\left(-cl\right) \leq \frac{1}{2},
\end{eqnarray}
thus
\begin{eqnarray}
||u^{\varepsilon}-u_{h,l}^{ms}||_{V}\leq Ch||f||_V,
\end{eqnarray}
where $C$ is a positive constant independent of the mesh size $h$.\hfill $\blacksquare$

\section{Numerical examples}
\begin{figure}
\begin{center}
\includegraphics[width=0.49\textwidth]{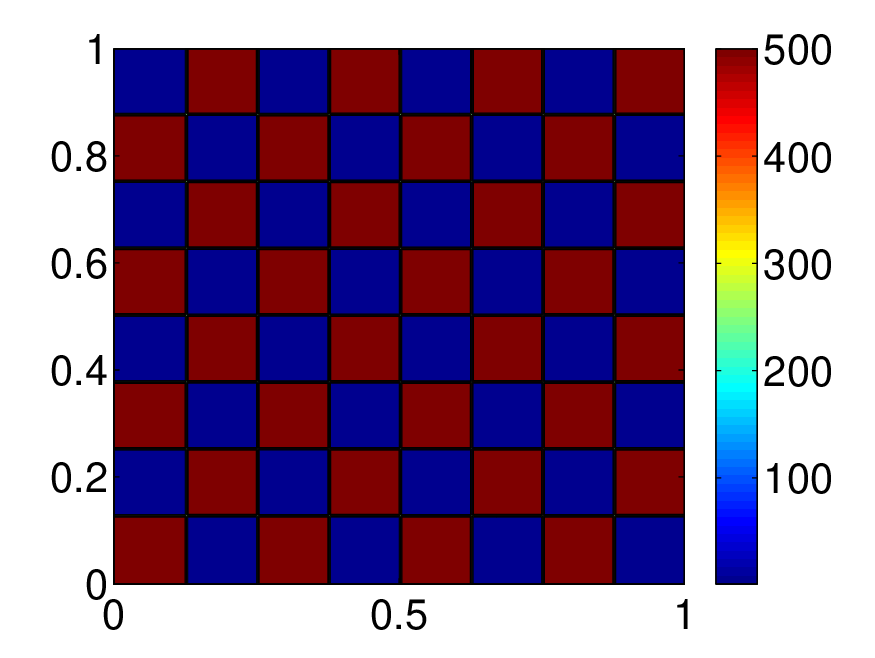}
\includegraphics[width=0.49\textwidth]{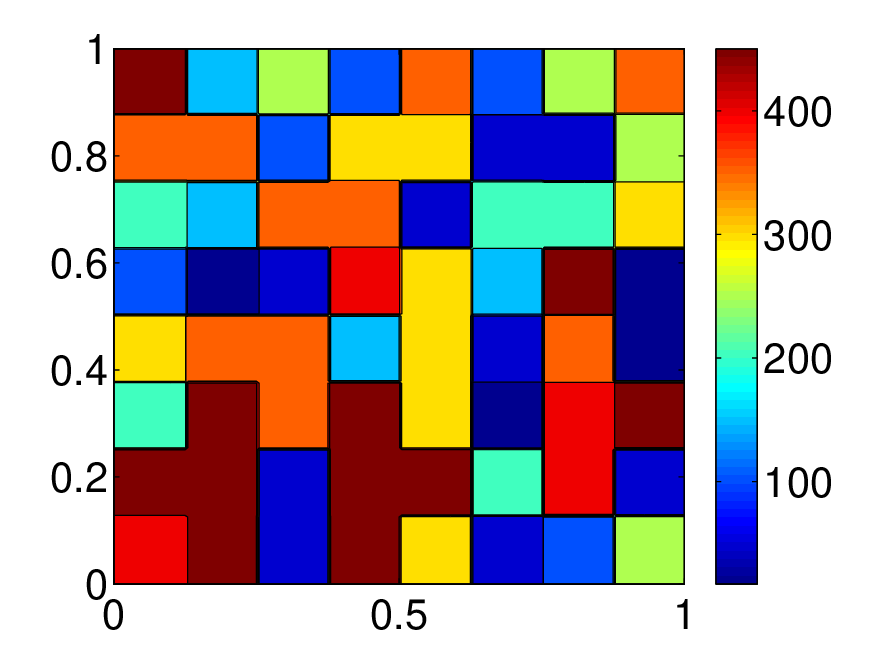}
\end{center}
\caption{Diffusion Coefficients: $A_1(x)$(left) and $A_2(x)$(right)}
\label{diffusion function}

\end{figure}

\begin{figure}
\begin{center}
\includegraphics[width=0.49\textwidth]{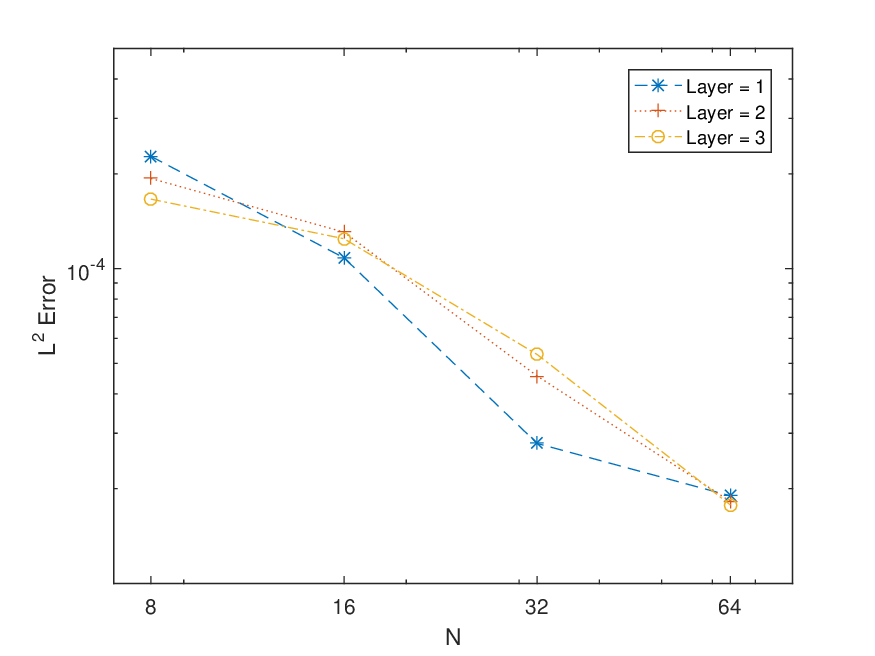}
\includegraphics[width=0.49\textwidth]{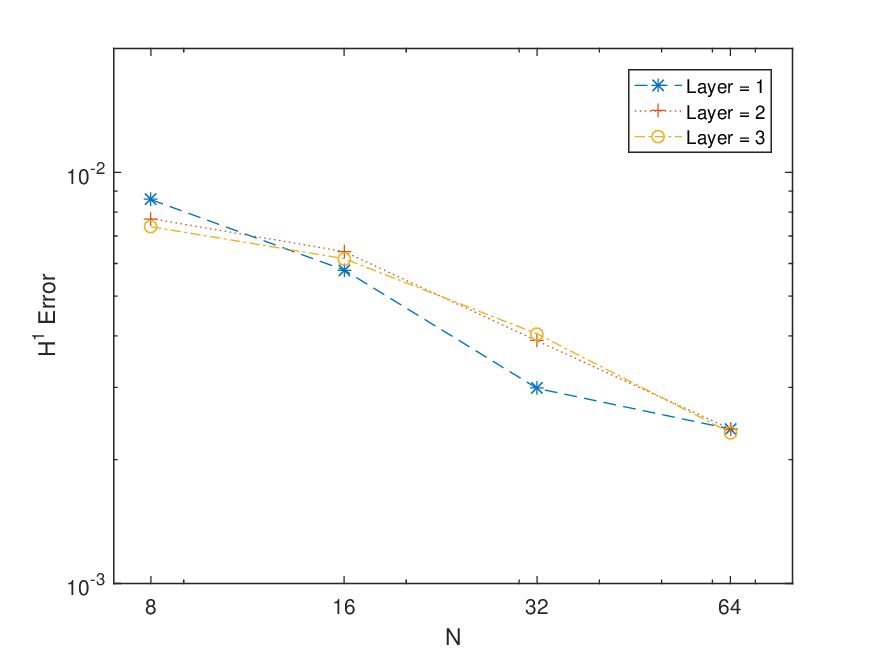}
\end{center}
\caption{$L^2$ error (left) and $H^1$ error (right) for diffusion coefficient $A_1(x)$.}
\label{erorr for A1}
\end{figure}

\begin{figure}
\begin{center}
\includegraphics[width=0.49\textwidth]{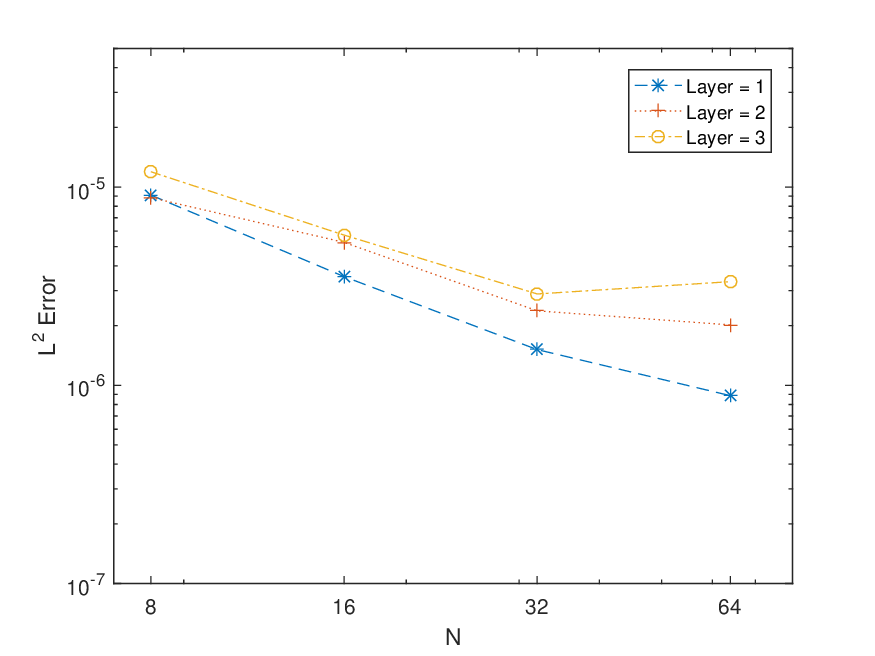}
\includegraphics[width=0.49\textwidth]{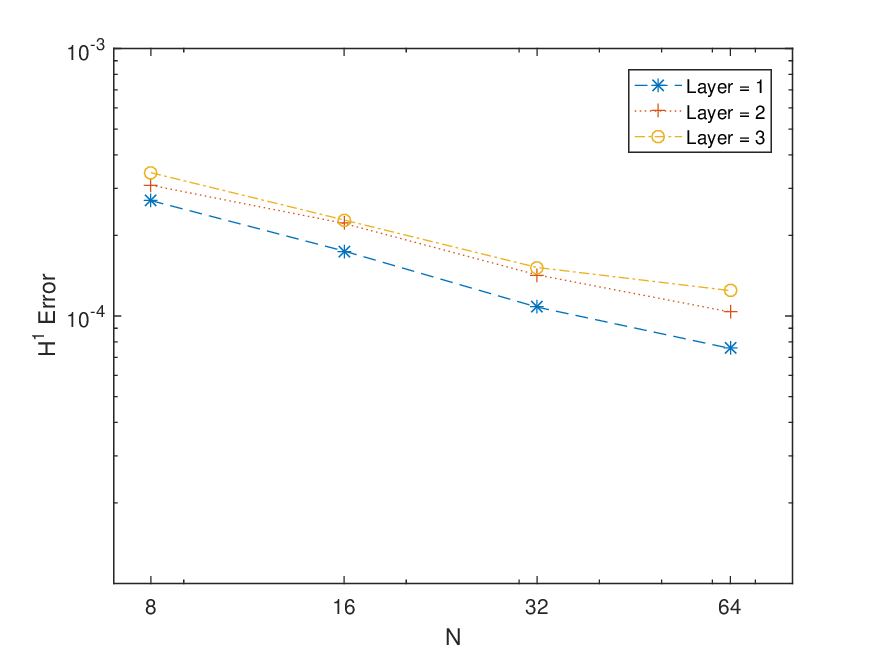}
\end{center}
\caption{$L^2$ error (left) and $H^1$ error (right) for diffusion coefficient $A_2(x)$.}
\label{erorr for A2}
\end{figure}

In this section, we propose some numerical experiments to verify the theoretical results. The domain is a unite square $\Omega=[0,1]\times[0,1]$, which is discretized with a uniform triangulation. Consider homogeneous boundary condition $u|_{\partial\Omega}=0$ and the outer force $f(x,y) = x$. The method is tested on two different types of diffusion coefficients, see Figure \ref{diffusion function}. The coarse mesh size is $H=1/8,1/16,1/32,1/64$ and the reference mesh $\mathcal{T}_h$ has width $h=1/256$. Since no analytical solutions are available, the standard finite element approximation $u_h$ on the reference mesh $\mathcal{T}_h$ serves as the reference solution. All fine scale computations are performed on subsets of $\mathcal{T}_h$. The number of layers of element patch is $l=1,2,3$. For diffusion coefficients $A_1$ and $A_2$, the $L_2$ errors and $H_1$ errors are displayed in Figure \ref{erorr for A1} and Figure \ref{erorr for A2}, respectively.


\section*{Acknowledgments}
This work was supported by the Natural Science Foundation of China (No. 62241203) and Science and Technology Research Project of Jiangxi Provincial Department of Education(No. GJJ211027).

%\input{bib}
%----------------------------------------------
%  This is the bib file for the paper vier.
%----------------------------------------------
%-------------------------REFERENCE-------------
\bibliographystyle{amsplain}

\end{sloppypar}
\end{document}